\documentclass[12pt,draft]{extartic}

\usepackage{amsmath,amsthm,amssymb,calc}
\usepackage[dvips]{graphics, graphicx}

\usepackage[english]{babel}

\def\eps{\varepsilon}

\newcounter{num}[section]

\newcommand{\Th}{\refstepcounter{num}
{\bf Theorem \arabic{section}.\arabic{num} }}

\newcommand{\Cor}{\refstepcounter{num}
{\bf Corollary \arabic{section}.\arabic{num} }}

\newcommand{\Note}{\refstepcounter{num}
{\it Note \arabic{section}.\arabic{num} }}

\def\eps{\varepsilon}
\def\_phi{\varphi}

\def\d{\delta}
\def\la{\lambda}
\def\ind{{\rm ind}}

\def\t{\tilde}
\def\o{\omega}
\def\ov{\overline}

\def\C{{\mathbb C}}

\def\r{\mathcal{R}}
\def\Z_N{{\mathbb Z}_N}
\def\Z{{\mathbb Z}}
\def\N{{\mathbb N}}

\def\l{\left}
\def\r{\right}

\author{Shkredov I.D.}
\title{ On monochromatic solutions of some nonlinear equations in $\Z/p\Z$
\footnote{ This work was supported
Pierre Deligne's grant based on his 2004 Balzan prize, grant RFFI
N 06-01-00383 and grant Leading Scientific Schools No.
691.2008.1}}
\date{}
\begin{document}
\maketitle

\begin{center}
 Annotation.
\end{center}

{\it \small
    Consider an arbitrary coloring  of integers with finite number of colors.
    Is it true that there are $x,y \in \Z$ such that $x+y$, $xy$ and $x$ have the same color?
    This is a well--known question of Ramsey theory has not solved yet.
    In the article we give a positive answer to the last question in the group
    $\Z/p\Z$, where $p$ is a prime number.
}
\\
\\
\\

\refstepcounter{section} \label{sec:introduction}

%\begin{center}
{\bf \arabic{section}. Introduction.}
%\end{center}
Let $k$ be a positive integer, and $\chi$ be an arbitrary coloring
of positive integers with $k$ colors. More formally, let $\chi :
\Z \to \{ 1,2,\dots, k \}$ be an arbitrary map and we associate
the segment of positive integers $\{ 1,2,\dots, k \}$ with $k$
different colors. If $f(x_1,\dots,x_n)=0$, $x_i \in \Z$ is an
equation, then a solution $(x^{(0)}_1,\dots,x^{(0)}_n)$ of the
equation is called {\it monochromatic} if all $x^{(0)}_i$ have the
same color. In other words, there exists $m\in \{1,\dots,k\}$ such
that $\chi(x^{(0)}_i) = m$, $i=1,\dots,n$. The problem of finding
monochromatic solutions (m.s. in short) of some equations was
considered in \cite{Ramsey_theory_book}--\cite{Deuber}.
%One of
The first result of the theory is Schur's theorem \cite{Schur}
which was proved in 1916.

\Th {\it
    For any coloring of $\Z$ there is a monochromatic solution of the equation
    $x+y=z$, where $x,y,z$ are non--zero numbers.
}

A classical result giving a complete answer about m.s. of any {\it
linear} equation is Rado's theorem (see \cite{Rado3}). M.s. of
linear equations and system of linear equations were studied not
only in $\Z$ but in another groups (see e.g. \cite{Deuber}). On
the other hand, we know a little about m.s. of  {\it nonlinear }
equations. For example, there is no answer yet to the question
about m.s. of the equation $x^2+y^2 = z^2$. We do not know are
there $x+y,xy$ having the same color for any\ finite coloring of
$\Z$ (see \cite{Imre_problems}, problem 3). It nothing to know
about more difficult problem on existing $x,y,x+y,xy \in \Z$ of
the same color. Note that Schur's theorem implies that there is
%existing of
a non--zero m.s. of the equation $xy=z$ (it is sufficient to
consider the set $\{ 2^n \}_{n\in \N_0}$).
% \bigsqcup \{ 2^{n \}_{n\in \N}$).
Concerning m.s. of non--linear equations including multiplication
see
%\cite{Bergelson_col1,Bergelson_col2,Bergelson_col3,non-linear1,non-linear2,non-linear3}.
\cite{Bergelson_col1}---\cite{non-linear3}. Unfortunately, in the
papers we have deal with some equations including either
multiplication or addition. In the article we positively answer to
the question from \cite{Imre_problems} in the case of the group
$\Z_p = \Z/p\Z$, where $p$ is a prime number. The last question
was asked to the author by Mathias Beiglb\"{o}ck.
%Mathias Beiglb\"{o}ck asked the author the last question.
Let us formulate our main result.

\Th {\it
    Let $p$ be a prime number, and $A_1,A_2,A_3 \subseteq \Z_p$ be any sets,
    $|A_1| |A_2| |A_3| \ge 40 p^{5/2}$.
    Then there are $x,y\in \Z_p$ such that $x+y\in A_1$, $xy \in A_2$, $x\in A_3$.
} \label{t:main_colouring_new}

%%?? ????? ????,
%????????? Theorem
%%????????? ??????????? ???????? ?? ?????? ????????? ?????????? --- ??????? \ref{t:main_colouring+} ????.

As one can see in the group $\Z_p$ we have not only a coloring
result but a statement with depends only on densities of our sets
$A_1$, $A_2$ and $A_3$ (so--called {\it density} result). In the
next section we give two intermediate variants of the theorem
above. Also we discuss some examples of analogous problems with
inclusions in the group $\Z_p$.

Theorem \ref{t:main_colouring_new} can be reformulate as follows.

\Cor {\it
    Let $p$ be a prime number, and $A_1,A_2,A_3 \subseteq \Z_p$ be any sets,
    $|A_1| |A_2| |A_3| \ge 40 p^{5/2}$.
    Then the equation
    \begin{equation}\label{f:three}
        a^2_1 + a_1 a_2 + a_3 = 0 \,,
    \end{equation}
    where $a_1\in A_1$, $a_2\in A_2$ and  $a_3\in A_3$ is always has a solution.
} \label{cor:main_colouring_new_reformulation}

In papers
%\cite{Shparlinskii_eq}, \cite{Sharkozy_eq}, \cite{Garaev_eq}
\cite{Shparlinskii_eq}---\cite{Garaev_eq} some equations with four
variables belonging to any sets $A_i$ of large cardinality (e.g.
$|A_i| \ge p^{1-\eps_0}$, where $\eps_0>0$ is an absolute
constant) were considered. Probably, equation (\ref{f:three}) is
the first example of equation with {\it three} variables which is
solvable for any sets  of cardinalities $\gg p^{1-1/6}$.

In paper \cite{Bourgain_more} (see also \cite{Ext_Bourgain_ext})
J. Bourgain studied a binary operation $\circ$ on $\Z_p$ which is
defined by the formula $x\circ y := x^2 + xy$. It was showed that
the operation has an "expansion"\, property. More precisely, if
$A,B \subseteq \Z_p$ are any sets $|A| \sim |B| < p^{1-\eps}$,
where $\eps>0$ is a constant, then
\begin{equation}\label{f:Bourgain_circ}
    |A \circ B | > |A|^{1+\d} \,,
\end{equation}
where $\d=\d(\eps) > 0$ is another constant. The value of
$\d(\eps)$ for $|A|\sim |B|$, $|A|,|B| \ge p^{\kappa}$, $\kappa >
1/2$ was calculated in \cite{Ext_Bourgain_ext}. Using the
operation $\circ$ we can reformulate our main theorem as follows
(see the proof in the next section).

\Cor {\it
    Let $p$ be a prime number, and $A,B \subseteq \Z_p$ be any nonempty sets.
    %%%, $|B| \ge \sqrt{p}$.
    Then
    $$
%        | A\circ B | \ge p - \frac{6 p^{7/4}}{(|A||B|)^{1/2}}  \,.
         | A\circ B | \ge (p - 1) - \frac{40 p^{5/2}}{|A||B|}  \,.
    $$
} \label{cor:new_Bourgain}

Thus Corollary \ref{cor:new_Bourgain} is an improvement of
(\ref{f:Bourgain_circ}),
%provided by
if $|A|,|B| \gg p^{3/4}$. A similar result was obtained by M.Z.
Garaev in \cite{Garaev_A(A+1)} for the operation $x*y := xy + x$.

The author is grateful to Mathias Beiglb\"{o}ck, Alexander Fish
for useful discussions and S.V. Konyagin for a number of helpful
advices and remarks. Also I acknowledge the Mathematical Sciences
Research Institute for its hospitality and providing me with
excellent working conditions.

\refstepcounter{section} \label{sec:proof}

{\bf \arabic{section}. The proof of the main result.} Let $f$ be
an arbitrary function from $\Z^*_p = (\Z / p\Z) \setminus \{ 0 \}$
to $\C$. By $\t{f} (z)$ denote the Fourier transform of $f$
\begin{equation}\label{f:Fourier_calc}
    \t{f} (z) = \sum_{x\in \Z^*_p} f(x) \ov{\chi_{z} (x)} \,.
\end{equation}
Here $\chi_z (x) = e^{2 \pi i z\cdot \ind x / (p-1)}$ is a
Dirichlet character. We put $\chi_z (0) = 0$ if $z\neq 0$. We will
use some well--known formulas
\begin{equation}\label{F_Par}
    \sum_{x\in \Z^*_p} |f(x)|^2 = \frac{1}{p-1} \sum_{\xi \in \Z^*_p} |\t{f} (\xi)|^2 \,.
\end{equation}
\begin{equation}\label{F_Par_sc}
    \sum_{x\in \Z^*_p} f(x) \overline{g(x)} = \frac{1}{p-1} \sum_{\xi \in \Z^*_p} \t{f}(\xi) \overline{\t{g}(\xi)} \,.
\end{equation}
\begin{equation}\label{f:inverse}
    f(x) = \frac{1}{p-1} \sum_{\xi \in \Z^*_p} \t{f}(\xi) \chi_{\xi} (x) \,.
\end{equation}
Let $A\subseteq \Z^*_p$ be a set. It is very convenient to write
$A(x)$ for the characteristic function of our set $A$. Thus
$A(x)=1$ if $x\in A$ and $A(x)=0$ otherwise.

First of all let us prove the following intermediate result.

\Th {\it
    Let $p$ be a prime number, and $A_1,A_2 \subseteq \Z_p$ be any sets,
    $|A_1| |A_2| \ge 20 p$.
    Then there are $x,y\in \Z_p$ such that $x+y\in A_1$, $xy \in A_2$.
} \label{t:main_colouring}

%Basically,
Theorem \ref{t:main_colouring} is an easy
%application
consequence of Weyl's bound for exponential sums with
multiplicative characters. Also, we prove a small modification of
Theorem \ref{t:main_colouring} in the section.

\Th {\it
    Let $p$ be a prime number, and $A_1, A_2 \subseteq \Z_p$  be any sets,
    $R\subseteq \Z^*_p$ be a multiplicative subgroup,
    and     $|R|^2 |A_1| |A_2| \ge 25 p^3 \log^2 p$.
    Then there are $x,y\in \Z_p$ such that $x\in R$, $x+y\in A_1$, $xy \in A_2$,
    and also there are $x,y\in \Z_p$ such that $x\in R$, $y-x\in A_1$, $xy \in A_2$.
} \label{t:main_colouring+}

\Note If we put $A_1=R$, $A_2=aR$, where $a$ is any non--zero
number, then the theorem above asserts that there are $x,y\in
\Z_p$ such that $x\in R$, $y-x\in R$, $xy \in aR$. Hence $R$ is a
basis
%in $\Z^*_p$
of order two (in other words $R+R=\Z^*_p$), provided by $|R| \ge 3
p^{3/4} \log^{1/2} p$. The last statement shows that our result is
sharp more or less.
%has some accuracy.
Indeed it is well--known (see e.g. \cite{KS}),                        %%%{\bf ???????????}
that a multiplicative subgroup $R$ is a basis of order two,
provided by $|R| > p^{3/4}$. This is the best result at the
moment. That is why one cannot replace the quantity $25 p^3 \log^2
p$ in Theorem \ref{t:main_colouring+} by, say, $o(p^3)$ (at least
using the current methods only).
%technic.
Anyway, one cannot replace $25 p^3 \log^2 p$ by $p^{2-\eps}$,
where $\eps>0$ is any number.

We will prove Theorems \ref{t:main_colouring},
\ref{t:main_colouring+} simultaneously. Let $\sigma$ be a number
of $x,y\in \Z^*_p$ such that $\o x+y\in A_1$, $xy \in A_2$, where
$\omega = \pm 1$. Let also $f(x) = \sum_y A_1(\omega x+y)
A_2(xy)$. Clearly, $\sigma = \sum_x f(x) = \t{f} (0)$. Let us find
Fourier coefficients of $f(x)$.
%Then
Using (\ref{f:inverse}), we have
$$
    \t{f} (z)
        =
            \sum_{x,y\in \Z^*_p} A_1(\omega x+y) A_2(xy) \ov{\chi_{z} (x)}
        =
            \sum_{x,\la\in \Z^*_p} A_1(x(\omega+\la)) A_2(\la x^2) \ov{\chi_{z} (x)}
                =
$$
$$
                =
                    \frac{1}{(p-1)^2} \sum_{r_1,r_2} \t{A}_1 (r_1) \t{A}_2 (r_2)
                        \sum_{x,\la} \chi_{r_1} ((\omega + \la)x) \chi_{r_2} (\la x^2) \ov{\chi_{z} (x)} \,,
$$
where $\t{A}_1 (r)$, $\t{A}_2 (r)$ are Fourier coefficients of the
sets $A_1,A_2$ (see formula (\ref{f:Fourier_calc})).
%Let $z\neq 0$.
Let $A'_1 = A_1 \bigcap \Z^*_p$, $A'_2 = A_2 \bigcap \Z^*_p$.
Obviously,  $|A'_1| \ge |A_1| - 1$, $|A'_2| \ge |A_2| - 1$. Using
the orthogonality property of the characters, we get
$$
    \t{f} (z) = \frac{1}{(p-1)} \sum_r \t{A}_1 (z-2r) \t{A}_2 (r) \sum_{\la} \chi_{z-2r} (\la+\omega) \chi_{r} (\la)
        =
$$
\begin{equation}\label{26.10.2008_1}
        =
            \frac{\t{A}_1 (z) |A'_2|}{(p-1)} \sum_{\la} \chi_{z} (\la+\omega)
                +
                    \frac{1}{(p-1)} \sum_{r\neq p-1}
                        \t{A}_1 (z-2r) \t{A}_2 (r) \sum_{\la} \chi_{z-2r} (\la+\omega) \chi_{r} (\la)\,.
\end{equation}
To estimate  the second term in (\ref{26.10.2008_1}) we use a
theorem from \cite{Johnsen}.

\Th {\it
    Let $m$ be a positive integer, $p$ be a prime number,
    $\chi_1, \dots, \chi_m$ be multiplicative characters,
    and $b_1,\dots, b_m \in \Z_p$ be different numbers.
    Then
\begin{equation}\label{est:Johnsen}
    \l| \sum_x \chi_1 (x+b_1) \chi_2 (x+b_2) \dots \chi_m (x+b_m) \r|
        \le
            (m - m_0 + 1) \sqrt{p} + m_0 + 1 \,,
\end{equation}
where $m_0$ is the number of the principal characters among
$\chi_1, \dots, \chi_m$. } \label{t:Johnsen}

Since $p\ge 2$, it follows that $3\sqrt{p}+1 \ge 2 \sqrt{p}+2$. We
have $|A_1| |A_2| \ge 20 p$ and $A_1,A_2 \subseteq \Z_p$. Hence
$p\ge 23$ and $|A_1|, |A_2| \ge 20$. Using the last inequality,
Theorem \ref{t:Johnsen}, Cauchy--Schwartz and Parseval's identity
(\ref{F_Par}), we obtain
$$
    \t{f} (z)
        =
            \frac{\t{A}_1 (z) |A'_2|}{(p-1)} \sum_{\la} \chi_{z} (\la+\omega)
                    + \theta_1 \frac{3\sqrt{p}+1}{(p-1)} p (|A_1| |A_2|)^{1/2}
                =
$$
\begin{equation}\label{f:t_f_Fourier}
                =
                    \frac{\t{A}_1 (z) |A'_2|}{(p-1)} \sum_{\la} \chi_{z} (\la+\omega) + \frac{7}{2} \cdot \theta_2 \sqrt{p} (|A_1| |A_2|)^{1/2} \,,
\end{equation}
where $|\theta_1|, |\theta_2| \le 1$. In particular, if we take $z
\equiv 0 \pmod {p-1}$ in formula (\ref{f:t_f_Fourier}) and use the
bound $|A_1| |A_2| \ge 20 p$, we get
\begin{equation}\label{f:av_f}
    \sigma = \t{f} (0) = |A_1| |A_2| + \t{\theta} \sqrt{19 p |A_1| |A_2| } > 0 \,,
\end{equation}
where $|\t{\theta}| \le 1$. This completes the prove of Theorem
\ref{t:main_colouring}.

Now let us obtain Theorem \ref{t:main_colouring+}. Note, that for
all $z\neq 0 \pmod {p-1}$, we have $|\sum_{\la} \chi_{z}
(\la+\omega)| \le 1$. It follows that for any such $z$, we get
\begin{equation}\label{f:zapas_future}
    |\t{f} (z)|
        \le
            \frac{|\t{A}_1 (z)| |A'_2|}{(p-1)} + 4 \sqrt{p} (|A_1| |A_2|)^{1/2}
                \le
                    5 \sqrt{p} (|A_1| |A_2|)^{1/2} \,.
\end{equation}
Consider the quantity
$$
    \sigma' = \sum_x R(x) f(x) = \frac{1}{p-1} \sum_z \ov{\t{R} (z)} \t{f} (z) \,.
$$
We have used (\ref{F_Par_sc}) to obtain the last formula. By
(\ref{f:t_f_Fourier}), Cauchy--Schwartz and a well--known bound
for Fourier coefficients of an arithmetic progression $\sum_z
|\t{R} (z)| \le p \log p$ (see e.g. \cite{Vinogradov_book}), we
have
$$
    \sigma'
        =
            \frac{|R||A'_1||A'_2|}{p-1} + \frac{\theta_3 |A'_2|}{(p-1)^2} \sum_z |\t{R} (z)| |\t{A}_1 (z)|
                +
                    \frac{4 \theta_2 \sqrt{p} (|A_1| |A_2|)^{1/2}}{p-1} \sum_z |\t{R} (z)|
                        =
$$
$$
    =
        \frac{|R||A'_1||A'_2|}{p-1} + \frac{\theta'_3 |A_2|}{(p-1)} (|R| |A_1|)^{1/2} +
            \frac{4 \theta_2 p\sqrt{p} (|A_1| |A_2|)^{1/2}}{p-1} \log p\,,
$$
where $|\theta_3|, |\theta'_3|, |\theta'_2| \le 1$. Using the
condition $|R|^2 |A_1| |A_2| \ge 25 p^3 \log^2 p$, we get $\sigma'
> 0$. This completes the prove of Theorem \ref{t:main_colouring+}.

\Note Careful analysis of the proof of the last theorem shows that
the inclusions  $x\in A$, $x+y\in A_1$, $xy\in A_2$ have place if
any of the sets $A$, $A_1$ or $A_2$ is a subgroup of $R$.

\Note One can generalize Theorem \ref{t:main_colouring}
%can be easily generalized
in the following way. For simplicity, let $A_1=A_2=A$. If $|A
\bigcap \Z^*_p|=\d p$, then simple average arguments shows that
there are $n\in \Z_p$ and $\la\in \Z^*_p$ such that $|A\bigcap
(A+n)| \ge \d^2 p$ and $|A\bigcap \la A| \ge \d^2 p$. Now we can
%Let us now
apply Theorem \ref{t:main_colouring} to these dense sets $A\bigcap
(A+n) \subseteq A$ and $A\bigcap \la A \subseteq A$. Clearly, we
can iterate the last procedure. For example,
%one can prove that
we get $x+y,xy,(x+y)z,xyz, x+y+w,xy+w,(x+y)z+w,xyz+w \in A$,
provided by the cardinality of $A$ is at least $\gg p^{7/8}$.
Obviously, there are very many similar examples.

Let us prove our main result. We use the method of successive
squaring in the spirit of papers \cite{Gow_4,Gow_m}.

{\bf Proof of Theorem \ref{t:main_colouring_new} and Corollary
\ref{cor:new_Bourgain}} Let $f(x)$ be the same function as above.
Formula (\ref{f:av_f}) gives us the mean value of the function.
Let us find the second moment of $f$. We have
$$
    \sigma_2 = \sum_x f^2 (x) = \sum_x |\sum_y A_1 (\omega x+y) A_2 (xy)|^2
        =
            \sum_{x\in \Z_p^*} |\sum_y A_2 (y) A_1 (\omega x+y x^{-1} )|^2 + A_2(0) |A_1|^2
                =
$$
\begin{equation}\label{f:second_m_f_1}
                =
                \sum_{y_1,y_2 \in A_2} \sum_{x\in  \Z^*_p} A_1 (\omega x+y_1 x^{-1}) A_1 (\omega x+y_2 x^{-1}) + A_2(0) |A_1|^2
                    =
                        \sigma'_2 + \sigma''_2 \,.
\end{equation}
Clearly, $|\sigma''_2| \le |A_1|^2$. Let us find the quantity
$\sigma'_2$. Let $\_phi (y_1,y_2) = \sum_{x\in \Z^*_p} A_1 (\omega
x+y_1 x^{-1}) A_1 (\omega x+y_2 x^{-1})$. Then
\begin{equation}\label{f:first_phi}
    \sum_{y_1,y_2} \_phi (y_1,y_2) = \sum_{x\in \Z^*_p} \sum_{y_1,y_2} A_1 (\omega x+y_1 x^{-1}) A_1 (\omega x+y_2 x^{-1})
        =
            (p-1) |A_1|^2 \,.
\end{equation}
Further
$$
    \sum_{y_1,y_2} \_phi^2 (y_1,y_2)
        =
            \sum_{y_1,y_2} \sum_{x_1,x_2 \in \Z_p^*}
                A_1 (\omega x_1+y_1 x_1^{-1}) A_1 (\omega x_1+y_2 x_1^{-1})
                A_1 (\omega x_2+y_1 x_2^{-1}) A_1 (\omega x_2+y_2 x_2^{-1}) \,.
$$
To calculate
%find
the second moment of the function $\_phi$ we need to find the
number of quadruples $(x_1,x_2,y_1,y_2)$ such that the product of
the characteristic functions of the set $A_1$ in the last formula
equals $1$. By $\mathcal{X}$ denote the set of such quadruples. In
other words, we can correspond  to any quadruple from
$\mathcal{X}$ a quadruple $(a_1,a_2,a_3,a_4)$, where $a_i \in
A_1$. We have
\begin{displaymath}
\left\{ \begin{array}{ll}
\omega x^2_1+y_1 = a_1 x_1 \\
\omega x^2_1+y_2 = a_2 x_1 \\
\omega x^2_2+y_1 = a_3 x_2 \\
\omega x^2_2+y_2 = a_4 x_2 \\
\end{array} \right.
\end{displaymath}
If $a_1 = a_2$, then $y_1=y_2$ and $a_3=a_4$. If we choose $y_1$,
$a_1 \in A_1$, $a_3 \in A_1$ in an arbitrary way, then solving
some quadratic equations we find $x_1$, $x_2$ in at least two
ways. It gives us at most $4 |A_1|^2 p$ solutions. In the case
$a_3=a_4$, we have $a_1=a_2$ and, consequently, all these
solutions were included in the situation above.
%$a_1=a_2$.
Now if $a_1=a_3$, then either $x_1=x_2$ or $x_1+x_2 = \omega a_1$
and, consequently, $y_1= \omega x_1 x_2$. Let us consider the
first situation. The equality $x_1=x_2$ implies $a_2=a_4$. It
follows that choosing $x_1$ and $a_1$, $a_2$ in an arbitrary way
we get $y_1$, $y_2$ unambiguously. It gives at most $|A_1|^2 p$
solutions. Let now $a_1=a_3$ and $x_1+x_2 = \omega a_1$.
%$y_1= \omega x_1 x_2$.
Then $y_1= \omega x_1 x_2$ and $x_1 (a_1-a_2) = x_2 (a_1-a_4)$.
Choosing $x_1$ and  $a_1$, $a_4$ in an arbitrary way we get
unambiguously $x_2$, $y_1$, $a_2$ and, as a consequence, $y_2$. It
gives at most $|A_1|^2 p$ solutions. The case $a_2=a_4$ adds  at
most $|A_1|^2 p$ solutions because the situation $x_1=x_2$ was
considered above. Clearly, the number of corresponding  quadruples
from $\mathcal{X}$ is at most $4 |A_1|^2 p + |A_1|^2 p + |A_1|^2 p
+ |A_1|^2 p = 7 |A_1|^2 p$. We can assume that $a_1\neq a_2$, $a_1
\neq a_3$, $a_3\neq a_4$, $a_2 \neq a_4$ and $y_1 \neq y_2$, $x_1
\neq x_2$. Subtracting the second equality from the first one and
the forth from the third, we have $x_1 = (y_1-y_2) / (a_1-a_2)$,
$x_2 = (y_1-y_2) / (a_3-a_4)$ and, as a consequence,
%consequently,
$x_1/x_2 = (a_3 - a_4) / (a_1-a_2)$. Further, we get
$$
    x_1 ( (a_3-a_4)^2 - (a_1-a_2)^2 ) = \omega (a_3-a_4) (a_2 a_3 - a_1 a_4)
$$
and
$$
    x_2 ( (a_3-a_4)^2 - (a_1-a_2)^2 ) = \omega (a_1-a_2) (a_2 a_3 - a_1 a_4) \,.
$$
Thus if $(a_3-a_4)^2 - (a_1-a_2)^2 \neq 0$, then we can find
$x_1$, $x_2$, and, as a consequence, $y_1$, $y_2$ by
$a_1,a_2,a_3,a_4$. If $(a_3-a_4)^2 - (a_1-a_2)^2 = 0$, then either
$a_3 - a_4 = a_1 - a_2$ or $a_3 - a_4 = -(a_1 - a_2)$. The first
alternative was considered before because we get  $x_1=x_2$. If
the second case holds then $x_1=-x_2$. But in the situation we
have $a_1=-a_3$ and choosing $x_1$, $a_1$, $a_2$ in an arbitrary
way we get $x_2$, $a_3$, $a_2$ and $y_1$, $y_2$ unambiguously. It
gives at most $p |A_1|^2$ quadruples from $\mathcal{X}$. Hence
\begin{equation}\label{f:second_phi}
    \sum_{y_1,y_2} \_phi^2 (y_1,y_2)
        \le
            |A_1|^4 + 8 p |A_1|^2 \,.
\end{equation}
Combining the last formula with (\ref{f:first_phi}), we obtain
\begin{equation}\label{tmp:30.09.2009_1}
    \sum_{y_1,y_2} \l( \_phi (y_1,y_2) - \frac{|A_1|^2}{p} \r)^2
        \le
            8 p |A_1|^2 + \frac{2|A_1|^4}{p}
                \le
                    10 p |A_1|^2 \,.
\end{equation}
Returning to (\ref{f:second_m_f_1}), we get
$$
    \sigma'_2
        =
            \sum_{y_1,y_2 \in A_2} \_phi(y_1,y_2)
                =
                    \frac{|A_1|^2 |A_2|^2}{p} + \sum_{y_1,y_2 \in A_2} \l( \_phi(y_1,y_2) - \frac{|A_1|^2}{p} \r) \,.
$$
Using Cauchy--Schwartz to estimate the second term of the last
formula, we have
$$
    \sigma'_2
        \le
            \frac{|A_1|^2 |A_2|^2}{p} + |A_1| |A_2| \sqrt{10 p} \,.
$$
Whence
$$
    \sigma_2 = \sum_x f^2 (x) \le \frac{|A_1|^2 |A_2|^2}{p} + |A_1| |A_2| \sqrt{10 p} + A_2 (0) |A_1|^2
$$
and we find the bound for the second moment of the function $f$.
Finally, consider the quantity
\begin{equation}\label{f:30.09.2009_2}
    \t{\sigma} = \sum_{x} A_3(x) f(x) = \frac{|A_1| |A_2| |A_3|}{p} + \sum_{x} A_3(x) \l( f(x) - \frac{|A_1| |A_2|}{p} \r) \,.
\end{equation}
By formula  (\ref{f:av_f}) and the estimate for $\sigma_2$, we
obtain
\begin{equation}\label{}
    \t{\sigma} \ge \frac{|A_1| |A_2| |A_3|}{p} - |A_3|^{1/2} \l( 2 |A_1|^{1/2} |A_2|^{1/2} p^{1/4} + A_2 (0) |A_1| + 3 |A_1|^{3/4} |A_2|^{3/4} p^{-1/4} \r)
        \ge
\end{equation}
\begin{equation}\label{f:30.09.2009_3}
        \ge
            \frac{|A_1| |A_2| |A_3|}{p} - |A_3|^{1/2} \l( 5 |A_1|^{1/2} |A_2|^{1/2} p^{1/4} + A_2 (0) |A_1| \r) \,.
\end{equation}
Clearly, we have the same upper bound for $\t{\sigma}$. Using
$|A_1| |A_2| |A_3| \ge 40 p^{5/2}$, we finally get
$$
    \t{\sigma} \ge \frac{|A_1| |A_2| |A_3|}{p} - 6 |A_1|^{1/2} |A_2|^{1/2} |A_3|^{1/2}  p^{1/4}
        > 0
        %\,.
$$
as required. To obtain  Corollary \ref{cor:new_Bourgain} one can
put $A_2$ equals the complement in $\Z_p \setminus \{ 0 \}$ to
$A\circ B$ and also put
%%%$A_1=-B$, $A_3=A$.
$\omega=-1$, $A_1=B$, $A_3=A$. After that
%one can
apply formulas (\ref{f:30.09.2009_2}) --- (\ref{f:30.09.2009_3}).
This completes the prove of the theorem and the corollary.

\Note There is a hope that using our method one can prove
solvability of the equation $p(a_1,a_2) = a_3$, where $a_i \in
A_i$, $A_i$, $i\in [3]$ are arbitrary sets of sufficiently large
density (see Corollary \ref{cor:main_colouring_new_reformulation})
and the polynomial $p(x,y)$ has the form $p(x,y) = p_1(x) + p_2
(x) g(y)$, where $g,p_1,p_2$ are some polynomials from $\Z_p [x]$,
$\deg g,p_1,p_2 \neq 0$ and $p_1, p_2$ are not affinely equivalent
(may be with some additional restrictions on $p_1,p_2$). It is
interesting to have  an example of a polynomial $p(x,y)$ of
another form such that the equation $p(a_1,a_2) = a_3$ is
solvable.

Let us note in conclusion that
%there is an open
a question remains open. Is it true that for any coloring of
$\Z^*_p$ exist $x,y,x+y,xy$ of the same color? Obviously, there is
no density result in the last question. To see this one can take
$A$ equals quadratic non--residuals or, as another example, the
set $\{ x\in \Z_p^* ~:~ x=2k+1,\, 0\le k\le (p-3)/4 \}$.

\end{document}